\documentclass[draft,10pt,oneside,intlimits]{amsart}
\usepackage[latin2]{inputenc}
\usepackage{amsmath}
\usepackage{amsthm}
\usepackage{amssymb}
\usepackage[mathscr]{eucal}
\usepackage{enumerate}

\def\NN{\mathbb{N}}
\def\RR{\mathbb{R}}

\def\smiley{\mbox{\hbox{\large$\bigcirc$}\hskip-0.95em\lower0.5pt\hbox{\tiny$\smallsmile$}\hskip-0.68em\raise4pt\hbox{..}}}

\DeclareMathOperator*{\diam}{diam}
\DeclareMathOperator*{\dd}{{d\hskip-1pt}}

\newcommand*{\bR}{{\mathbb R }}
\def\RRb{\overline{\mathbb R}}
\newcommand*{\bC}{{\mathbb C }}

\newcommand*{\fM}{{\mathfrak{M}} }
\newcommand*{\fMe}{\ensuremath{ {\mathfrak{M}}_1  }}

\newcommand*{\ve}{\varepsilon}

\newcommand*{\ws}{weak$^*$}

\def\conv{\mathop{\mathrm{conv}}}

\def\cconv{\mathop{\mathrm{\overline{conv}}}}

\theoremstyle{plain}
\newtheorem{proposition}{Proposition}
\newtheorem{theorem}[proposition]{Theorem}
\newtheorem{lemma}[proposition]{Lemma}
\newtheorem{corollary}[proposition]{Corollary}

\theoremstyle{definition}
\newtheorem{definition}[proposition]{Definition}

\theoremstyle{definition}

\newtheorem{remark}[proposition]{Remark}
\newtheorem{question}[proposition]{Question}

\DeclareMathOperator*{\supp}{supp}

\newcounter{rem}
\newcounter{rev}
\begin{document}
\numberwithin{proposition}{section}
\numberwithin{equation}{section}


\title[Rendezvous numbers of metric spaces]
{Rendezvous numbers of metric spaces \\ -- a potential theoretic
approach}

\author[B.~Farkas]{B\'{a}lint Farkas}
\address{Technische Universit\"{a}t Darmstadt\newline Fachbereich Mathematik, AG4 \newline Schlo\ss gartenstra\ss
e 7\newline D-64289 Darmstadt, Germany }
\email{farkas@mathematik.tu-darmstadt.de}
\author[Sz.~Gy.~R\'{e}v\'{e}sz]{Szil\'{a}rd Gy. R\'{e}v\'{e}sz$^*$}
\address{Alfr\'{e}d R\'{e}nyi Institute\newline
Hungarian Academy of Sciences\newline Re\'{a}ltanoda u.~13--15\newline
H-1053, Budapest, Hungary} \email{revesz@renyi.hu}


\thanks{
The present publication was supported by the Hungarian-French
Scientific and Technological Governmental Cooperation, project \#
F-10/04 and the Hungarian-Spanish Scientific and Technological
Governmental Cooperation, project \# E-38/04.}

\begin{abstract}
The present work draws on the understanding how notions of
general potential theory -- as set up, e.g., by Fuglede -- explain
existence and some basic results on the ``magical'' rendezvous
numbers. We aim at a fairly general description of rendezvous
numbers in a metric space by using systematically the potential
theoretic approach.

In particular, we generalize and explain results on invariant
measures, hypermetric spaces and maximal energy measures, when
showing how more general proofs can be found to them.
\end{abstract}

\maketitle \let\oldfootnote\thefootnote
\def\thefootnote{*}
\footnotetext{Supported in part by Hungarian National Foundation
for Scientific Research, Grant \# T 049301.}
\let\thefootnote\oldfootnote
\section{Introduction}\label{sec:Intro}
It was proved by O.~Gross that for a compact, connected metric
space $(X,d)$ there exists a unique number $r=r(X)$ such that for
every finite point system $x_1,\dots,x_n\in X$, $n\in \NN$ one
always finds an $x\in X$ with
\begin{equation}\label{eq:AV}
r=\frac{1}{n}\sum_{i=1}^n d(x,x_n)~.
\end{equation}
Such a number is called the \emph{rendezvous number} of the
metric space $X$. Since the first result of Gross \cite{Gross},
rendezvous numbers have been attracting much attention and been
generalized in many directions: considering weak rendezvous
numbers (Thomassen \cite{Th}), replacing the metric by some
continuous symmetric function (Stadje \cite{Stadje}) or
considering instead of the finite average in \eqref{eq:AV} the
mean value with respect to some probability measure $\mu$ (Elton,
Cleary, Morris, Yost \cite{CMY}). In such abstract investigations
various minimax principles play important role. (See also Morris,
Nickolas \cite{MN}, Nickolas, Yost \cite{NY}, Stranzen
\cite{Str}).

Our aim is to put the investigations on the existence and
uniqueness of rendezvous numbers in the framework of abstract
potential theory, which has been around since the 60s, but
apparently has not gained its due recognition in this field. In
this paper we continue \cite{FR} with the study of related notions
such as invariant measures and maximal energy (Wolf \cite{W3},
Bj\"{o}rck \cite{Bj}). It turns out that general principles such as
existence of capacitary measures or Frostman's equilibrium
theorem are accounted for the existence of invariant measures. In
the past ten years or so, generalizations to possibly infinite
dimensional, hence not locally compact spaces, appeared. In
particular, calculation of rendezvous numbers of unit spheres in
Banach spaces fascinated many authors. (See, e.g., Baronti,
Casini, Papini \cite{baronti/casini/papini:2000},
Garc\'{\i}a-V\'{a}zquez, Villa \cite{GVV}, Lin \cite{Lin}, Wolf
\cite{W1}).

\vskip2ex First we spend some words on technicalities and recall
the appropriate setting of potential theory in locally compact
spaces. For convenience we add $+\infty$ to the set of real
numbers, i.e., let $\overline{\bR}:=\bR\cup\{+\infty\}$ endowed
with its natural topology such that $\overline{\bR}_+$ will be
compact. Moreover, we will use the notation $\conv E$ for the
convex hull of a subset $E\subset \RR$ and $\cconv E$ for the
closed convex hull in $\overline{\RR}_+$, meaning, for example,
$\overline{\conv}(0,+\infty)=[0,+\infty]$.

In the fundamental work of Fuglede \cite{Fug}, general potential theory is presented in locally
compact spaces. So unless otherwise stated $X$ will be a locally compact Hausdorff space.
Nevertheless we would like to use the developed tools also on metric spaces. To this end the
appropriate results will be carried over to this case in Section \ref{sec:metric}.

On the space $X$ we will consider, a usually fixed, kernel function $k$ in the sense of Fuglede
\cite[p.~149]{Fug}. That is, we assume that $k:X\times X\rightarrow \RRb$  is {\em lower
semicontinuous} (l.s.c.) as a two variable function over $X\times X$, and that $-\infty < k(x,y)
\le +\infty$. Moreover, we assume that $k\ge 0$, and that $k$ is symmetric, i.e., $k(x,y)=k(y,x)$
for all $x,y\in X$.

Denote by $\mathfrak{M}(X)$ the family of positive, finite,
regular Borel measures on $X$; $\fMe(X)$ will stand for the
subset of probability measures. For any $H\subset X$ we let
$$\fMe(H):=\{\mu \in \fMe(X)
:\:\mbox{$H$ is $\mu$-measurable and $\mu(H)=\mu(X)$}\}~.
$$
The customary topology on $\fM$ is the vague topology which is the locally
convex topology determined by the seminorms $\mu\mapsto \left|\int_X f \dd\mu\right|$, $f\in
C_c(X)$. In most places we consider only the family $\fMe(K)$ of probability measures supported
on the same compact set $K$. In this case, by the Riesz Representation Theorem,
$\fM(K)$ is the positive cone of $C(K)'$, and the weak$^*$-topology determined by $C(K)$ and the vague topology
coincide.
\subsection{Potential and energy}
Just as in the classical case, the \emph{potential} and
\emph{energy} of $\mu\in \fM$ are
\[
    U^\mu(x) := \int_X  k(x,y)\dd\mu(y) ~, \qquad
    W({\mu}) := \int_X\int_X  k(x,y)\dd\mu(y)\dd\mu(x)~.
\]

\begin{definition}\label{def:qq}
Let $H\subset X$ be fixed, and $\mu\in\fMe(X)$ be arbitrary. First put
\begin{equation}\label{d:QmuH}
Q(\mu,H):=\sup_{x\in H} U^{\mu}(x)~,\qquad\text{and also}\qquad
\underline{Q}(\mu,H):=\inf_{x\in H} U^{\mu}(x)~.
\end{equation}
For any two sets $H,L\subset X$ we define the quantities
\begin{align}\label{qstarcapacity}
q(H,L)&:= \inf_{\mu\in\fMe(H)} Q(\mu,L) &\text{and}&& \underline{q}(H,L)&:=\sup_{\nu\in\fMe(H)}
\underline{Q}(\nu,L)~.
\end{align}
\end{definition}

\begin{definition} For $\mu\in\fMe$, recalling
Fuglede \cite{Fug}, we write
\begin{align}\label{mucapacities}
W(\mu):&= \sup_{K\Subset X} W({\mu_K})~,
\end{align}
with $\mu_K$ denoting the trace of $\mu$ on the set $K$. The
Wiener energy (reciprocal capacity) of any set $H\subset X$ is
\begin{align*}
w(H):=&\inf_{\mu\in\fMe(H)} W(\mu)=\inf_{\mbox{\scriptsize$\mu\in\fMe(H)$}\atop \mbox{\scriptsize$\supp \mu \Subset H$}} W(\mu) =
    \inf_{K\Subset H} w(K)~, 
\end{align*}
where for the last forms see \cite[(2), p.~152]{Fug}.
\end{definition}
We remind that Fuglede \cite[(1), p.~152]{Fug}
defines the so-called ``uniform'' and ``de la Vall\'{e}e-Poussin''
energies $U(\mu):=Q(\mu,X)$ and $V(\mu):=Q(\mu,\supp\mu)$ and
their counterparts $u(H)$ and $v(H)$ for subsets of $H\subset X$.
In \cite{Bela} their relationship to the Chebyshev constant (see
below) and transfinite diameter is studied. However, we will not
need these special cases, in the following.

We will use the following statement from \cite[Lemma 2.2.1]{Fug}
or \cite[Lemma 1.5]{FR}.
\begin{lemma}\label{lem:lsc}Let $H,L\subset X$. The functions
below are lower semicontinuous
\begin{enumerate}[a)]
\item \hfil$\displaystyle\fM_1(H)\times L\ni(\mu,x)\mapsto U^{\mu}(x)\Bigl(:=\int_X
k(x,y)\dd\mu(y)\Bigr)~,$\hfill
\item \hfil$\displaystyle\fMe(H)\ni\mu\mapsto Q(\mu;L)\Bigl(:=\sup_{x\in L}
U^{\mu}(x)\Bigr).~$\hfill
\end{enumerate}
\end{lemma}

\subsection{Chebyshev contants}

\begin{definition} For arbitrary $H,L \subset X$ we define the
\emph{(general) $n^{\text{th}}$ Chebyshev constant} of $L$ with
respect to $H$ as
\[
M_n(H,L):= \sup_{w_1,\ldots,w_n \in H} \inf_{x\in L}
    {\frac{1}{n}}\Bigg( \sum_{k=1}^{n} k(x,w_k) \Bigg),
\]
and the \emph{(general) $n^{\text{th}}$  dual Chebyshev constant}
of $L$ relative to $H$ as
\[
    \overline{M}_n(H,L):= \inf_{w_1,\ldots,w_n \in H} \sup_{x\in L}
    {\frac{1}{n}}\Bigg(
   \sum_{j=1}^{n} k(x,w_j)
   \Bigg)
   ~.
\]
The \emph{$n^{\text{th}}$
Chebyshev constant} of $H$ is  $M_n(H):=M_n(H,H)$ and the
\emph{$n^{\text{th}}$ dual Chebyshev constant} of $H$ is  $\overline{M}_n(H) :=
\overline{M}_n(H,H)$.
\end{definition}

The following proposition is proved by showing that the
respective sequences are quasi-mo\-no\-to\-nous, hence Fekete's
lemma (see \cite{Fek}, or also \cite{Bela}, \cite{FR}) applies.
\begin{proposition}
For any $H,L\subset X$, the Chebyshev constants $M_n(H,L)$ and $\overline{M}_n(H,L)$ converge, more precisely
\begin{equation*}
\sup_{n\in\NN} M_n(H,L)=\lim_{n\rightarrow \infty} M_n(H,L)\qquad\mbox{and}\qquad\inf_{n\in\NN}\overline{M}_n(H,L)=\lim_{n\rightarrow \infty} \overline{M}_n(H,L).
\end{equation*}
\end{proposition}

The limits whose existence is assured by the previous proposition
are denoted by $M(H,L)$ and $\overline{M}(H,L)$ (and for $H=L$
also by $M(H)$, $\overline{M}(H)$), respectively.

\subsection{Rendezvous intervals}
We define the (weak) rendezvous number(s) and average distance number(s) of the
space $X$, or even of subsets of $X$.
Again, for good reasons (explained in more detail in
\cite{FR}) we define these notions in dependence of two sets as
variables.

\begin{definition}\label{def:randi}For arbitrary subsets $H, L \subset X$ the $n^{\text{th}}$
\emph{(weak) rendezvous set} of $L$ with respect to $H$ is
\begin{equation}\label{nthrandi}
R_n(H,L):=\bigcap_{w_1,\dots,w_n\in H} \cconv
\Bigl\{p_n(x):=\frac 1n \sum_{j=1}^n k(x,w_j)~~:~~ x\in L
\Bigr\},\quad R_n(H):=R_n(H,H)~.
\end{equation}
Correspondingly, one defines
\begin{align}\label{randidef}
R(H,L)&:=\bigcap_{n=1}^\infty R_n(H,L)~,& R(H):=R(H,H). \intertext{Similarly, one defines the
\emph{(weak) average set} of $L$ with respect to $H$ as} \label{avridef}
A(H,L)&:=\bigcap_{\mu\in\fMe(H)} \cconv \Bigl\{U^{\mu}(x)~~:~~ x\in L \Bigr\} ~,& A(H):=A(H,H).
\end{align}
\end{definition}

\begin{remark}\label{AmH} Denoting the interval
\begin{equation}\label{Amuint}
A(\mu,L):=[\underline{Q}(\mu,L), Q(\mu,L)]= \cconv\{U^{\mu}(x)~:~x\in L\}~,
\end{equation}
we see that $R_n(H,L)$, $R(H,L)$ and $A(H,L)$ are all of the form
$\bigcap_{\mu} A(\mu,H)$, with $\mu$ ranging over all averages of
$n$ Dirac measures at points of $H$, over all measures finitely
supported in $H$ and having only rational probabilities, and over
all of $\fMe(H)$, respectively, see \cite{FR}.

\end{remark}
\begin{remark}\label{sdc}
If $k$ is a continuous kernel -- in particular when it is a
metric -- and $L$ is compact, then it suffices to take $\conv$
instead of $\cconv$, since then together with $k$ also
$U^{\mu}(x)$ is continuous for any $\mu$. Thus for compact
subsets $L$ of metric spaces a real number $r\in\RR_{+}$ belongs
to $R(H,L)$ if and only if for any $x_1,\dots,x_n \in H$
($n\in\NN$) we always have points $y,z\in L$ satisfying
\begin{equation}\label{nthweakrendezvous}
\frac 1n \sum_{j=1}^n k(y,x_j)\le r \qquad\text{and}\qquad \frac 1n \sum_{j=1}^n k(z,x_j) \ge
r~,
\end{equation}
which is the usual definition of \emph{weak rendezvous numbers} in metric spaces.

Moreover, in case the set $L$ is connected, this is further
equivalent to the existence of a ``rendezvous point'' $x\in L$
with
\begin{equation}\label{nthstrongrendezvous}
\frac 1n \sum_{j=1}^n k(x,x_j) = r~.
\end{equation}
In particular, for compact and connected $L$ in a metric space
(or in a locally compact space with continuous kernel $k$) the
rendezvous set $R(H,L)$ consists of a \emph{unique} point, say
$R(H,L)=\{r(H,L)\}$, if this latter property
\eqref{nthstrongrendezvous} is satisfied only for $r=r(H,L)$.
\end{remark}
\begin{remark}\label{sd} If $k$ is only l.s.c., also potentials are l.s.c.,
which entails that they take on their infimum over compact sets. Thus for compact $L$ the first
half of the above equivalent formulation \eqref{nthweakrendezvous} remains valid even for
general kernels. However, for the second part we must already write that ``$\forall s<r ~
\exists z\in L$ such that $\frac 1n \sum_{j=1}^n k(z,x_j) > s$''. Such modification of the
formulation is necessary also when we consider sets $L\subset X$ which are not compact, or when
we are discussing the case when $+\infty\in R(H,L)$. Clearly, in our settings $R_n(H,L)$,
$R(H,L)$ and $A(H,L)$ are subsets of $[0,\infty]$, but note that traditionally rendezvous
numbers or average numbers are considered only among the reals.
\end{remark}

With the above notions at hand, the following description of
various rendezvous sets is easy to see, cf.~\cite{FR}.
\begin{proposition}\label{prop:reni}For arbitrary subsets $H, L \subset X$
we have
\begin{align}\label{randneqchebyn}
&R_n(H,L)= [M_n(H,L),\overline{M}_n(H,L)]~, && R_n(H)= [M_n(H),\overline{M}_n(H)]~,\notag \\
&R(H,L)= [M(H,L), \overline{M}(H,L)]~,  && R(H)= [M(H), \overline{M}(H)]~, \\
&A(H,L)=[\underline{q}(H,L),q(H,L)]~,  &&A(H)=[\underline{q}(H),q(H)]~ .\notag
\end{align}
\end{proposition}

\begin{remark}\label{rem:convent}
Note that intervals appearing in proposition \ref{prop:reni}
\emph{may indeed be empty}, meaning, for example, that  $q(H,L)<\underline{q}(H,L)$, cf.~\cite{FR} and also Theorem \ref{th:existence} below.
\end{remark}

\section{General results on rendezvous numbers}\label{sec:generalities}

The following theorem, known as Frostman's theorem in the
classical case, shows the relationship between the potential of a
capacitary (energy-min\-imi\-zing) measure and the Wiener energy
of a set. See \cite[Theorem 2.4]{Fug}, or \cite{FR}.
\begin{theorem}[\bf Fuglede\rm]\label{frostmantetele} Let $k$ be a
positive, symmetric kernel and $K \Subset X$ be a compact set with
$w(K)<+\infty$. Every $\mu\in \fMe(K)$ having minimal energy
($W(\mu)=w(K)$) satisfies
\begin{align}
\label{FrostA}
U^\mu(x) & \ge w(K) \mbox{ for nearly every } x \in K~, \\
\label{FrostB}
U^\mu(x) & \le w(K) \mbox{ for every } x \in \supp \mu~, \\
\label{FrostC} U^\mu(x) & = w(K) \mbox{ for } \mu\mbox{-almost every } x \in X~.
\end{align}
\end{theorem}

\begin{remark}\label{noexception} In case $k$ is continuous, or even if
only it is bounded on $K\times K$, there can be no sets of finite
measure but infinite energy. Therefore, the exceptional set in
\eqref{FrostA} (which refers to probability measures of $\fMe
(K)$) must be void, and \eqref{FrostA} \emph{holds everywhere}.
\end{remark}

The following results are recalled from \cite{FR}.

\begin{theorem}\label{cor:sum} Let $H, L\subset X$, then
\begin{align}\label{eq:sum1}
 M(H,L)\leq \underline{q}(H,L)\leq&\: q(L,H)\leq \overline{M}(L,H)~.\\
 \intertext{If $L\subset X$ is compact, then}
M(H,L)= \underline{q}(H,L)=&\:q(L,H)~. \label{eq:sum2} \\
\intertext{If $K\subset X$ is compact and $k$ is continuous, then}
 &\:q(L,K)=\overline{M}(L,K)~. \label{eq:sum3}
 \end{align}
\end{theorem}

\begin{theorem}\label{th:existence} Let $X$ be a locally compact Hausdorff
space, $\emptyset \ne H \subset L \subset X$ be arbitrary, and let $k$ be any nonnegative,
symmetric kernel on $X$. Then the intervals $R_n(H,L)$, $R(H,L)$ and $A(H,L)$ are nonempty.
\end{theorem}

\begin{remark} In general, $A(H)\subsetneqq R(H)$ is possible, see \cite[Remark 6.4]{FR}.
Also, the rendezvous
intervals can be ``almost empty'': consider, e.g.,
$R_n(\RR,\RR)=\{+\infty\}$. This and Remarks \ref{sdc} and
\ref{sd} already explain the slightly disturbing situation that
some papers state that ``there is no rendezvous number'' for cases
where we find one. However, not only $+\infty$ can show up in the
closure of intervals for the definition of rendezvous numbers,
hence not only $+\infty$ can be a rendezvous number for us while
does not exist for other authors. See \cite{FRnorm} for the cases of
$\ell_p$ spaces.
\end{remark}

\begin{theorem}\label{th:unique} Let $X$ be any locally compact Hausdorff
topological space, $k$ be any l.s.c., nonnegative, symmetric kernel function, and $K\Subset X$ be
compact. Then $A(K)$ consists of one single point. Furthermore, if $k$ is continuous and $K$ is
compact, then even $R(K)$ consists of only one point.
\end{theorem}

\begin{theorem}\label{thm:AeqR}
If $k$ is continuous and $L$ is compact, we have $R(H,L)=A(H,L)$ for all $H\subset X$.
\end{theorem}

In general, the theory of rendezvous numbers seems to be
flourishing in the context of metric spaces instead of locally
compact spaces with a Fuglede-type kernel. The latter theory is
more general regarding the kernel, but is a bit restrictive in
requiring local compactness of the underlying space. This gap is
filled by indicating that the above potential theoretical
approach works even for metric spaces, even if not locally
compact. That leads us to the next section.

\section{Rendezvous numbers for metric spaces}\label{sec:metric}

Note that for a nonnegative, Borel measurable (e.g., a continuous
or l.s.c.) function $f:X\to \RR$ and a (positive, finite) Borel
measure $\mu$ the integral $\int_X f\dd\mu$ may be defined as a
Lebesgue integral. Thus the potential $U^\mu(x):=\int_X d(x,y)\dd
\mu(y)$ -- and hence all related notions, considered previously --
are defined (cf.~Section \ref{sec:Intro}). Further, keeping the
notations from Section \ref{sec:Intro}, we have that
$\mu\in\fMe(H)$ implies that $H$ is $\mu$-measurable.

These remarks are already sufficient to define the Chebyshev
constants, rendezvous intervals and to show the equalities
\eqref{randneqchebyn} as well as Theorem \ref{cor:sum} and
\ref{th:existence} in the metric space setting. We will further
elaborate on this matter in \cite{FRnorm} regarding normed spaces.
Now we cover the theory of rendezvous numbers of metric spaces.

\label{ss:Elton}

Gross' result on the existence of a rendezvous number was
generalized by G.~Elton to general Borel probability measures in
place of finite convex combinations of Dirac measures. Note that
for continuous kernels on compact sets the closure can be skipped
from Definition \ref{def:randi} (as in that case potentials are
continuous, and a continuous image of a compact set is always
closed). On the other hand, Thomassen \cite{Th} extended the
result to not necessarily connected spaces by considering
so-called \emph{weak rendezvous numbers}, which is equivalent to
applying the convex hull in the definition, cf.~Remark \ref{sdc}.
Hence in our notation merging Gross', Thomassen's and Elton's
theorems corresponds to the following result.

\begin{theorem}[\bf Gross--Thomassen--Elton\rm]\label{th:metric}
Let $(X,d)$ be a compact metric space. Then we have
$A(X)=R(X)=\{r(X)\}$. Furthermore, there exist probability
measures $\mu,\nu \in \fMe$ with the property that
\begin{equation}\label{Elton}
U^{\mu}(x) \leq r(X) \le U^{\nu}(y) \qquad (\forall x, y \in X)~.
\end{equation}
\end{theorem}

\begin{remark} As mentioned a couple of times above, by
compactness and continuity here we have exactly the same result
even if closure is skipped from Definition \ref{def:randi};
furthermore, if the space $X$ is connected, then neither is any
need for convex hull.
\end{remark}

A further extension is due to Stadje \cite{Stadje}, who
essentially obtained the assertion of Theorem \ref{th:unique} concerning $R(X)$. He in fact assumed
connectedness, but this assumption is easily removed when
considering weak rendezvous numbers, i.e., convex hulls of values
attained by the respective potential functions; also, he
considered only $R(X)$, and not $A(X)$. We see that all
these results follow from Theorem \ref{th:unique}.

Note that Elton did not publish his result, but references to his
work \cite{CMY, MN} mention that he proved his statement even for
continuous, nonnegative and symmetric functions $f$ (in place of
the metric $d$) over compact connected Hausdorff topological
spaces. In any case, his results are now included in the
following.

\begin{theorem} Let $X$ be a locally compact Hausdorff space,
$k$ a symmetric, l.s.c., nonnegative kernel, and $\emptyset\neq
H\subset X$ be arbitrary, while $\emptyset\neq K\Subset X$ be
compact subsets of $X$. Then there exists $\mu \in \fMe(K)$ with
the property that
\begin{equation}\label{Eltonup}
U^{\mu}(x) \leq q(K,H) \qquad (\forall x \in X)~
\end{equation}
and for all $\epsilon>0$ there exists $\nu \in \fMe(K)$ with
\begin{equation}\label{Eltondowneps}
\underline{q}(K,H)-\ve \leq U^{\nu}(y) \qquad
(\forall y \in X)~.
\end{equation}
Moreover, if the kernel $k$ is continuous and bounded on $K\times
H$, then we have
\begin{equation}\label{Eltondown}
\underline{q}(K,H) \leq U^{\nu}(y) \qquad (\forall y
\in X)~.
\end{equation}
\end{theorem}
\begin{proof} By definition, there exist
$\mu_n\in \fMe(K)$ with $Q(\mu_n,H) \leq q(K,H) + 1/n$. Since
$\fMe(K)$ is \ws-compact by compactness of $K$, there exists a
subnet $\mathcal{N}$ of these measures converging to some $\mu
\in \fMe(K)$. In view of lower semicontinuity (see Lemma
\ref{lem:lsc} b) ), $q(K,H) \geq \liminf_{\mathcal N} Q(\mu_n,H)
\geq Q(\mu,H)$, hence the assertion \eqref{Eltonup}.

Inequality \eqref{Eltondowneps} is just the definition.

To prove \eqref{Eltondown} consider the ``dual'' kernel $\ell:=
C-k$ whenever $k$ is continuous and bounded by some constant $C$.
Then $\ell$ is nonnegative, symmetric and l.s.c., and the first
part applies. It is easy to check that to any measure
$\nu\in\fMe(K)$ the potentials with respect to $k$ and $\ell$ are
related\footnote{We will use subscript notation for the kernel.
For instance, we write $U^\mu_k$ to emphasize that the potential
$U^\mu$ is understood with respect to the kernel $k$.} by
$U_k^{\nu}=C-U_{\ell}^{\nu}$, while
$\underline{q}_k(K,H)=C-q_{\ell}(K,H)$.
\end{proof}

Note that we did not assume $H$ to be compact. However, in case
we have a pair of compact sets $K,L$, then a continuous kernel is
necessarily bounded on $K\times L$ and thus \eqref{Eltondown}
follows. In particular, for $K=L$ and a continuous kernel Elton's
result is obtained using also $\underline{q}(K)=q(K)$, i.e., the
last part of Theorem \ref{th:unique}.

\section{Invariant measures and rendezvous
numbers}\label{ss:invariance}

Following Morris and Nickolas \cite{MN}, but extending the notion
from $H=L=X$ to arbitrary subsets $H,L\subset X$, and from
metrics $d$ to arbitrary kernels $k$, we call a measure
$\mu\in\fMe(H)$ \emph{$k$-invariant} (on $L$), if the respective
potential integral is constant:
\begin{equation}\label{kinvariant}
U^{\mu}_k(x):=\int_X k(x,y) \dd\mu(y) \equiv \text{const.}\qquad
\text{(for all $x\in L$)}~.
\end{equation}
Saying only that $\mu$ is \emph{$k$-invariant} refers to the
central case $L=X$. Then an extension of the result of Morris and
Nickolas \cite{MN} to general kernels $k$ sounds as follows.

\begin{theorem}\label{th:kinvari} Assume
that $X$ is a locally compact Hausdorff topological space and $k$
is a nonnegative, l.s.c., symmetric kernel function. Let $\emptyset\neq H\subset
L\subset X$ be arbitrary and assume that there exists a measure
$\mu\in\fMe(H)$ which is $k$-invariant on $L$. Then we have
$$
A(H,L)=A({\mu},L)~.
$$
Furthermore, if $k$ is continuous and $L$ is compact, then we
even have
$$
R(H,L)=A({\mu},L)~.
$$
\end{theorem}
\begin{proof} Note that $A({\mu},L)$ being the (convex
closure of the) set of values of $U^{\mu}(x)$ when $x$ runs $L$, invariance immediately implies
that $\# A({\mu},L)=1$. Hence only (non-empty) existence of \eqref{avridef} is needed to conclude
$A(H,L)=A({\mu},L)$: this follows from Theorem \ref{th:existence}.
To obtain the last assertion from this, it suffices to refer to Theorem \ref{thm:AeqR}.
\end{proof}

\begin{corollary}[\bf Morris--Nickolas\rm]\label{cor:dinvari}
Let $(X,d)$ be a compact (connected) metric space. Assume that
there exists a $d$-invariant measure $\mu_0\in\fMe(X)$. Then we
have
$$
A(X)=R(X)=\{r(X)\}=A({\mu_0},X), \qquad U^{\mu_0}(x)\equiv r(X)
\quad (\forall x\in X)~.
$$
\end{corollary}

\section{Maximal energy and rendezvous numbers}\label{sec:Wolf}

Wolf \cite{W3} presents a theory of rendezvous numbers and
maximal (i.e., maximal energy) measures on compact connected
metric spaces $(X,d)$. Let us revise these results in this
section. Following Bj\"{o}rck \cite{Bj}, Wolf says that a
probability measure $\mu_1\in \fMe(X)$ is \emph{maximal}, and
that the space has \emph{maximal energy} $E(X)$, if
\begin{equation}\label{maxenergy}
E(X):=E_d(X):=\sup_{\mu\in \fMe(X)} W_d({\mu})=W_d({\mu_1})~.
\end{equation}
By weak$^*$-compactness of $\fMe(X)$, existence of $\mu_1$ is
obvious. Wolf proves that $r(X)\leq E(X)$, and also gives examples when
$r(X)<E(X)$.

\begin{theorem}[\bf Wolf\rm]\label{th:Wolf} Let $(X,d)$ be a compact
metric space. Then
\begin{itemize}
\item[(i)] $r_d(X)\leq E_d(X)$.
\item[(ii)] If $r_d(X)= E_d(X)$, then there exists some
$d$-invariant measure in $\fMe(X)$.
\end{itemize}
\end{theorem}

In his proof in \cite[pp.~396--397]{W3} Wolf uses properties of metrics rather heavily. Here we
extend the result first proving the following.

\begin{theorem}\label{th:kinvariant} Let $\emptyset\neq K\subset L\subset X$ be arbitary sets and $k$ be a nonnegative, l.s.c.,
symmetric kernel. Then we have
\begin{equation}\label{RKgeqwK}
\min R(K,L) \geq w (K)~. \end{equation} In particular, if $k$ is
continuous and $K$ is compact, then the set $K$ has a unique
rendezvous number $r(K)$, and we have
\begin{equation}\label{rKwK}
r(K) \geq w (K).
\end{equation}
Furthermore, if $r(K)= w(K)$, then there exists some
$k$-invariant measure in $\fMe(K)$.
\end{theorem}
\begin{proof} Existence of rendezvous numbers
$A(K,L) \subset R(K,L)$ are provided by Theorem \ref{th:existence},
and we also know $R(K,L)=[M(K,L),\overline{M}(K,L)]$ (see Proposition
\ref{prop:reni}). At this point \eqref{RKgeqwK}
follows from the fact that $M(K)\ge w(K)$ and that $M(K,L)\geq M(K)$ (see \cite{Bela, FR}).

According to Theorem \ref{th:unique}, continuity of $k$ on the
compact set $K$ implies uniqueness of the rendezvous numbers
$A(K)=R(K)=\{r(K)\}$, giving the second part of the statement.

Furthermore, let now $r(K)= w(K)$ be assumed. Since $k$ is
continuous and $K$ is compact, in this case $w(K)<+\infty$ is
obvious.

Take now a probability measure $\mu\in \fMe(K)$ minimizing
$Q(\mu,K)$, i.e., with $Q(\mu,K)=q(K)=r(K)$. Such a measure
exists, because $\mu\mapsto Q(\mu,K)$ is l.s.c.~in view of Lemma
\ref{lem:lsc} b), and $\fMe(K)$ is weak$^*$-compact. For any such
$\mu$ we have
$$
w(K)\leq W(\mu)=\int_K U^\mu(x)\dd\mu(x)\leq \sup_{x\in K} U^\mu
(x)=Q(\mu,K)=r(K)=w(K),
$$
so equality must hold throughout. Hence $\mu$ minimizes also
$W(\mu)$ (it is a capacitary measure). For this $\mu$ the
inequality \eqref{FrostA} of Theorem \ref{frostmantetele} holds,
moreover, it holds \emph{everywhere} on $K$ by Remark
\ref{noexception}. But then $w(K)\leq U^\mu(x) \leq Q(\mu,K) =
r(K)=w(K)$, hence $U^\mu(x)=w(K)$ holds for all $x\in K$, and
$\mu$ is seen to be a $k$-invariant measure.
\end{proof}

Now we are in the position to deduce Wolf's theorem as an easy
corollary.

\begin{proof}[Proof of Wolf's Theorem \ref{th:Wolf}]
Let $k:=\diam(X) - d$, which is then a continuous, symmetric,
nonnegative kernel function. By the previous Theorem
\ref{th:kinvariant}, $r_k(X)\geq w_k(X)$. Also,
$E_d(X)=\diam(X)-w_k(X)$ is immediate. In view of Theorem
\ref{th:unique}, uniqueness of the rendezvous numbers hold both
with respect to $d$ and $k$, and thus we have
$r_k(X)=q_k(X)=\underline{q}_k(X)$ and also
$r_d(X)=q_d(X)=\underline{q}_d(X)$. Definition \ref{def:qq}
immediately yields $q_k(X)=\diam(X)-\underline{q}_d(X)$ and
$q_d(X)=\diam(X)-\underline{q}_k(X)$, which show (i).

Let us now assume $E_d(X)=r_d(X)$, i.e., $r_k(X)=w_k(X)$. As then a
$k$-invariant measure $\mu\in\fMe(X)$ exists, and obviously
$U^{\mu}_k(x)=\diam(X)-U^{\mu}_d(x)$, the very same measure is also
$d$-invariant and even (ii) follows.
\end{proof}

Wolf also treats the converse question: when does the existence
of a $d$-invariant measure imply the equality of the maximal
energy and the rendezvous number? He uses the following notion.

\begin{definition} A metric space $(X,d)$ is called \emph{hypermetric},
if for all finite collections of points $x_i\in X$ ($i=1,\dots,n$)
and real scalars $c_i$ ($i=1,\dots,n$) with $\sum_{i=1}^n c_i=0$,
we have
\begin{equation}\label{hypermetric}
\sum_{i=1}^n \sum_{j=1}^n c_i c_j d(x_i,x_j) \le 0~.
\end{equation}
\end{definition}

Wolf discusses how the notion proves to be useful, a number of
well-known spaces being hypermetric spaces; for the details see
\cite{W3}. Here we only aim at revealing the potential theoretic
background even of this notion and the corresponding converse
result of Wolf.

Observe that by the densness of convex linear combinations of Dirac
measures in $\fMe$ for the weak$^*$-topology (see, e.g.,
\cite[Proposition 2.1.2, page 52]{nieberg:1991} and the
Kre\u{\i}n--Milman Theorem) and in view of continuity of $d$,
\eqref{hypermetric} implies that we also have
\begin{equation}\label{morehyper}
W_d({\mu-\nu})= \int_X \int_X d(x,y)
\dd(\mu-\nu)(x)\dd(\mu-\nu)(y) \le 0 \qquad \left( \mu,\nu \in
\fMe \right) ~.
\end{equation}
Translating this property to a property of the ``dual kernel''
$k:=m-d$, where $m:=\max_{X\times X}d$ is constant (in fact, the
diameter), we get
\begin{equation}\label{khyper}
W_k({\mu-\nu})= \int_X \int_X k(x,y)
\dd(\mu-\nu)(x)\dd(\mu-\nu)(y) \ge 0\qquad \left( \mu,\nu \in
\fMe \right)~.
\end{equation}
This property is almost identical with the notion of (positive)
definiteness, having great importance in potential theory, see
\cite[p.~151]{Fug}. Fuglede calls a kernel $k$ (positive)
definite, if for any signed regular Borel measure $\sigma\in
\fM^{\pm}(X)$ one has $W_k({\sigma})\geq 0$. This is slightly more
stringent, than \eqref{khyper}, where only
$\|\sigma^+\|=\|\sigma^-\|$ is considered, but \eqref{khyper}
will suffice in the next argument.

\begin{theorem}\label{th:converses} Let $\emptyset\neq K\subset L\subset X$ be arbitrary sets.
Assume \eqref{khyper} and that $A(K,L)=\{a(K,L)\}$. If there is a
probability measure $\mu_0\in \fMe(K)$ which is $k$-invariant on
$L$, then we have $a(K,L)=w(K)$ and $U^{\mu_0}$ is constant
$w(K)$ (on $L$).
\end{theorem}
\begin{proof}
Since $A(K,L)=\{a(K,L)\}$, we have that $a(K,L)\in A(\mu,L)$ for all $\mu\in\fMe(K)$, so $U^{\mu_0}(x)=a(K,L)$ for all $x\in L$.

Applying \eqref{khyper} for $\nu:=\mu_0$ and arbitrary $\mu\in
\fMe(K)$, we obtain
\begin{align}\label{kposdef}
0 \leq W({\mu-\mu_0}) & = \int_K
\left(U^{\mu}(y)-a(K,L)\right)\dd(\mu-\mu_0)(y)=\int_K
U^{\mu}(y)\dd(\mu-\mu_0)(y)= \\ &=W({\mu})- \int_K U^{\mu_0}(x)
\dd\mu(x)=W({\mu})-a(K,L) \qquad \left( \mu\in \fMe \right)~, \notag
\end{align}
hence for all $\mu\in \fMe(K)$ we have $W({\mu}) \geq a(K,L)$.
Taking infimum over all $\mu\in \fMe(K)$ yields $w(K)\geq
a(K,L)$. On the other hand, \eqref{RKgeqwK} of Theorem
\ref{th:kinvariant} and  $R(K,L)\supseteq A(K,L)$  yield
$a(K,L)\geq w(K)$, hence $a(K,L)=w(K)$, and also $U^{\mu_0}\equiv
w(K)$.
\end{proof}

Using that for continuous kernels one always has the uniqueness of the rendezvous numbers and the equality $A(K,L)=R(K,L)$, we arrive to the following corollary.
\begin{corollary}\label{cor:convers} Let $\emptyset\neq K$ be a compact set and $k$ a continuous kernel.
Assume \eqref{khyper} (or its equivalent discrete form, analogous to \eqref{hypermetric}).
If there is a probability measure $\mu_0\in \fMe(K)$ which is $k$-invariant on
$K$, then we have $r(K)=w(K)$ and $U^{\mu_0}$ is constant $w(K)$ (on $K$).
\end{corollary}

\begin{corollary}[\bf Wolf\rm]\label{th:conversewolf} Let the compact
metric space $(X,d)$ be hypermetric. If there is a $d$-invariant
probability measure $\mu\in \fMe(X)$, then we have $r(X)=E(X)$.
Furthermore, the potential of the $d$-invariant measure is
constant $r(X)$ and is of maximal energy.
\end{corollary}
\begin{proof}
Note that for compact $X$, the inequalities  \eqref{hypermetric}
and \eqref{morehyper} are equivalent. Again we consider the
continuous, symmetric, nonnegative kernel function $k:=\diam(X) -
d$. By Theorem \ref{th:converses}, $r_k(X)=w_k(X)$, and so
$r_d(X)=E_d(X)$ follows. Moreover, $U^{\mu}_k(x)\equiv w_k(K)$
implies the rest of the statement.
\end{proof}
\begin{question} Does there exist a true invariant measure for,
e.g., the unit sphere $S_{\ell_p}$?
\end{question}
\begin{question} Do we have an Elton-type ``separation theorem''
even in not locally compact spaces? In normed spaces?
\end{question}
\begin{definition} A measure $\mu\in\fMe(H)$ is termed
``$\ve$-quasi-invariant on $L$'' if $\sup_L U^{\mu}-\inf_L
U^{\mu}\leq \ve$.
\end{definition}
\begin{question} If the rendezvous number is unique, do we
have a (quasi-) converse: Do there exist at least
$\ve$-quasi-invariant measures?
\end{question}
This is interesting as there is way to conclude the argument of Theorem \ref{th:converses} from
the very existence of such $\ve$-quasi-invariant measures.
\begin{proposition} Let $X$ be any (not necessarily locally compact)
Hausdorff topological space, and $H,L\subset X$ be arbitrary with
$A(H,L)\ne \emptyset$. Assume that for all $\ve>0$ there exists
some $\ve$-quasi-invariant measure on $L$ from $\fMe(H)$. Take any
sequence $\ve_n\to 0$ ($n\to \infty$) together with the
corresponding measures $\mu_n\in \fMe(H)$, $\ve_n$-quasi-invariant
on $L$, and consider any values $\rho_n$ attained by the
respective potentials $U^{\mu_n}$ on $L$. We then have $\rho_n
\to a(H,L)$ as $n\to \infty$, where the average number exists
uniquely, i.e., $A(H,L)=\{a(H,L)\}$
\end{proposition}
\begin{proof} By $\ve_n$-quasi-invariance, $A(\mu_n,L)\subset [\rho_n-\ve_n,
\rho_n+\ve_n]$. As the intersection of the sets $A(\mu_n,L)$
contains $A(H,L)$, the intersection must be nonempty by
condition. Therefore, the intersection is a diameter $0$ nonempty
subset -- that is, a single point -- of $\RR$. However, as this
set $\{\rho\}$ contains the nonempty set $A(H,L)$, we conclude
$\rho=a(H,L)$. It is clear that $\rho_n\to \rho$ as $n\to\infty$.
\end{proof}

\begin{remark}
The analogue of the above proposition for the rendezvous numbers
also hold, where $R(H,L)$, $r(H,L)$ replace $A(H,L)$ and $a(H,L)$
respectively.
\end{remark}

\parindent0pt

\end{document}